\documentclass[12pt,draft]{extartic}

\usepackage{amsmath,amsthm,amssymb,calc}
\usepackage[dvips]{graphics, graphicx}

\usepackage[english]{babel}

\newcounter{num}[section]

\newcommand{\Th}{\refstepcounter{num}
{\bf Theorem \arabic{section}.\arabic{num} }}

\newcommand{\Lemma}{\refstepcounter{num}
{\bf Lemma \arabic{section}.\arabic{num} }}

\newcommand{\Pred}{\refstepcounter{num}
{\bf Proposition \arabic{section}.\arabic{num} }}

\newcommand{\Cor}{\refstepcounter{num}
{\bf Corollary \arabic{section}.\arabic{num} }}

\newcommand{\Exm}{\refstepcounter{num}
{\bf Example \arabic{section}.\arabic{num} }}

\newcommand{\St}{\refstepcounter{num}
{\bf Statement \arabic{section}.\arabic{num} }}

\newcommand{\Def}{\refstepcounter{num}
{\it Definition \arabic{section}.\arabic{num} }}

\newcommand{\Proof}{{\bf Proof. }}

\def\eps{\varepsilon}
\def\_phi{\varphi}

\def\a{\alpha}
\def\d{\delta}
\def\l{\lambda}

\def\v{\vec}
\def\F{\widehat}
\def\L{\Lambda}
\def\m{\times}
\def\t{\tilde}
\def\o{\omega}
\def\ov{\overline}
\def\z{{\mathbb Z}}
\def\C{{\mathbb C}}
\def\r{\mathcal{R}}
\def\Z_N{{\mathbb Z}_N}
\def\Z{{\mathbb Z}}
\def\N{{\mathbb N}}

\author{Shkredov I.D.}
\title{On sets of large exponential sums, I
\footnote{This work was supported by RFFI grant no.
06-01-00383, President's of Russian Federation grant N 1726.2006.1
and INTAS (grant no. 03--51--5-70).}}
\date{}
\begin{document}
\maketitle

\begin{center}
 Annotation.
\end{center}

{\it \small Let $A$ be a subset of $\mathbb{Z} / N\mathbb{Z}$
and let $\mathcal{R}$ be the set of large Fourier coefficients of $A$.
Properties of $\mathcal{R}$ have been  studied in works of M.--\,C. Chang and
B. Green.
Our result is the following  :
the number of quadruples $(r_1, r_2, r_3 , r_4 ) \in \mathcal{R}^4$
such that
$r_1 + r _2  =  r_3 + r_4$  is at least $|\mathcal{R}|^{2+\epsilon}, \epsilon>0$.
This statement shows that the set $\mathcal{R}$ is highly structured.
We also  discuss some of the generalizations and applications
of our result.
}
%\end{center}
\\
\\
\\

\refstepcounter{section}

{\bf \arabic{section}. Introduction.}

Let $N$ be a positive integer.
By $\Z_N$ denote the set $\z / N \z$.
Let $f: \Z_N \to \C$ be an arbitrary function.
Denote by $\F{f}$ the Fourier transform of $f$
\begin{equation}\label{}
    \F{f} (r) = \sum_{n\in \Z_N} f(n) e(-nr) \,,
\end{equation}
where $e(x) = e^{-2\pi i x/N}$.
By Parseval's identity
\begin{equation}\label{f:Par}
    \sum_{r\in \Z_N} |\F{f} (r)|^2  = N \sum_{x\in \Z_N} |f(x)|^2 \,.
\end{equation}
 Let $\d,\a$ be real numbers, $0<\a \le \d \le 1$ and
 let $A$ be a subset of $\Z_N$ of cardinality $\d N$.
 It is very convenient to write $A(x)$ for such a function.
 Thus $A(x) = 1$ if $x\in A$ and $A(x) = 0$ otherwise.
 Consider the set $\r_\a$ of large exponential sums of the set $A$
\begin{equation}\label{f:R_def}
    \r_\a = \r_\a (A) = \{~ r\in \Z_N ~:~ |\F{A} (r)| \ge \a N ~\} \,.
\end{equation}
 In many problems of combinatorial number theory is important to know the structure of the set $\r_\a$
 (see \cite{Gow_surv}).
% Для  многих задач комбинаторной теории чисел важно знать структуру множества $\r_\a$.
 In other words what kind of properties $\r_\a$ has?

 It is easy to see that $0\in \r_\a$ and $\r_\a = -\r_\a$.
 Further, using (\ref{f:Par}) we obtain $|\r_\a| \le \d / \a^2$.
 Are there any non--trivial properties of the set $\r_\a$?

 In 2002 M.--\,C. Chang proved the following result \cite{Ch_Fr}.

 \Th {\bf (Chang)}
 \label{t:Chang}
 { \it
    Let $\d,\a$ be real numbers, $0<\a \le \d \le 1$,
    $A$ be a subset of $\Z_N$, $|A| = \d N$.
    Then there exists a set $\L = \{ \lambda_1, \dots, \lambda_{|\L|} \} \subseteq \Z_N$,
    $|\L| \le 2 (\d / \a)^2 \log ( 1/\d)$ such that
    for any $r \in \r_\a$ we have
    \begin{equation}\label{f:presentation}
        r = \sum_{i=1}^{|\L|} \eps_i \lambda_i  \pmod N \,,
    \end{equation}
    where
    $\eps_i \in \{ -1,0,1 \}$.
 }

 Using approach of paper \cite{Ruzsa_Freiman} (see also  \cite{Bilu})
 Chang applied her result to prove the famous Freiman's theorem \cite{Freiman}
 on sets with small doubling.
 Let us formulate this beautiful result.

 The set $Q \subseteq \Z$
 $$
    Q = \{ n_0 + n_1 \l_1 + \dots + n_d \l_d ~:~ 0 \le \l_i < m_i \}
 $$
 is said to be a {\it $d$--dimensional arithmetic progression}.

 \Th {\bf (Freiman)}
 \label{t:Freiman}
 {\it
    Let $C>0$ be a real number and  $A \subseteq \Z$
    be an arbitrary set.
    Suppose that $|A+A| \le C|A|$.
    Then there exist numbers $d$ and $K$ depend on $C$ only and
    a $d$--dimensional arithmetic progression $Q$ such that $|Q| \le K |A|$
    and $A\subseteq Q$.
 }

 The second application of Theorem \ref{t:Chang} was obtained by B. Green in paper \cite{GreenA+A}
 (see also \cite{BourgainA+A,FreimanA+A} and \cite{CrootA+A}).
 Let us formulate one of the main results of \cite{GreenA+A}.

%% \Th {\bf (Грин)}                                                                   % м.б. выкинуть
%% \label{t:Green_A+A}
%% {
%% }

 \Th {\bf (Green)}
 \label{t:Green_3A}
 { \it
    Let $A$ be an arbitrarily subset of $\Z_N$, $|A| = \d N$.
    Then $A+A+A$ contains an arithmetic progression of length at least
    \begin{equation}\label{f:Green_3A}
        2^{-24} \d^5 ( \log(1/\d) )^{-2} N^{\d^2 / 250 \log (1/\d)} \,.
    \end{equation}
 }

 %%%Об оценке сверху на длину максимальной арифметической прогрессии в суммах множеств
 %%%см. работу \cite{RuzsaA+A}.

 In paper \cite{Green_Chang_exact} Green showed that Chang's theorem is sharp in a certain sense.
 Let $E = \{ e_1, \dots, e_{|E|} \} \subseteq \Z_N$ be an arbitrary set.
 By ${\rm Span}(E)$ denote the set of all sums $\sum_{i=1}^{|E|} \eps_i e_i$, where
 $\eps_i \in \{ -1,0,1 \}$.

 \Th {\bf (Green)}
 \label{t:Green_Chang_exact}
 { \it
    Let $\d,\a$ be real numbers, $\d \le 1/8$, $0< \a \le \d/32$.
    Let also
    \begin{equation}\label{}
        \left( \frac{\d}{\a} \right)^2 \log (1/ \d) \le \frac{\log N}{\log \log N} \,.
    \end{equation}
    Then there exists a set $A\subseteq \Z_N$, $|A| = [ \d N ]$ such that
    the set $\r_\a$ does not contain in ${\rm Span}(\L)$
    for any set $\L$ of cardinality $2^{-12} (\d/ \a)^2 \log (1/ \d)$.
 }

 If the parameter $\a$ is close to $\d$ then the structural properties of the set $\r_\a$
 was studied in papers
 \cite{Freiman_Yudin,Besser,Lev}, see also survey \cite{Konyagin_Lev}.

 The main result of the paper is the following theorem.

 \Th
 \label{t:main}
 {\it
    Let $\d,\a$ be real numbers, $0< \a \le \d$,
    $A$ be a subset of $\Z_N$, $|A| = \d N$
    and
    $k\ge 2$
    Suppose that $B\subseteq \r_\a \setminus \{ 0 \}$ is an arbitrary set.
    Then the number
    \begin{equation}\label{f:T_k_def}
        T_k (B) := |\{ ~ (r_1,\dots, r_k, r_1', \dots, r_k') \in B^{2k} ~:~
                    r_1 + \dots + r_k = r_1' + \dots + r_k' ~ \}|
    \end{equation}
    is at least
    \begin{equation}\label{f:T_k}
        \frac{\d \a^{2k}}{2^{4k} \d^{2k}} |B|^{2k} \,.
    \end{equation}
 }

 Let us show that the statement of Theorem  \ref{t:main} is not trivial
 in the case when
 $\d$ tends to zero as $N$ tends to infinity
 (if $\d$ does not tends to zero as $N \to \infty$ then there are not non--trivial
 restrictions on structure of the set $\r_\a$,
 see papers \cite{Katznelson,Nazarov,Ball}).
 Let us consider the simplest case $k=2$.
 Let the cardinality of $\r_\a$ is equal to $\Theta(\d/\a^2)$.
 Using Theorem \ref{t:main}, we obtain that the number of solutions of the equation
 \begin{equation}\label{f:1}
    r_1 + r_2 = r_3 + r_4\,, \quad \mbox{ where } \quad r_1,r_2,r_3,r_4 \in \r_\a \setminus \{ 0 \}
 \end{equation}
 equals $\Theta(\d /\a^4)$.
 There are three sorts of trivial solutions of equation (\ref{f:1}).
 The first sort of solutions is $r_1 = r_3$, $r_2=r_4$, the second --- $r_1=r_4$, $r_2=r_3$, and
 the third --- $r_1 = -r_2$, $r_3 = -r_4$.
 So the number of trivial solutions of (\ref{f:1}) does not exceed $3|\r_\a|^2$.
 Since the cardinality of $\r_\a$ does not exceed $\d/\a^2$ it follows that
 $3|\r_\a|^2$ at most $3 \d^2 /\a^4$.
 This number is less than $\d / \a^4$ as $\d$ tends to zero.
 Thus Theorem \ref{t:main} asserts that equation (\ref{f:1})
 has non--trivial solutions.
 In the sense Theorem \ref{t:main} shows that the set $\r_\a$
 has some additive structure.

 We prove Theorem \ref{t:main} in section \ref{proofs}.
 In \S \ref{matrix} we obtain a matrix generalization of our approach.
 We use Gowers uniformity norms (see \cite{Gow_m}) in our proof.

 In section  \ref{applications} we obtain some applications of our main result.
 We show that Theorem and W. Rudin's inequality implies Chang's result.
 Moreover we derive an improvement of Theorem \ref{t:Chang}
 (see Theorem \ref{t:Chang_log}).
 Also we obtain an application of Theorem \ref{t:main} to Freiman's Theorem \ref{t:Freiman}.

 In our forthcoming papers we are going to obtain further
 results on sets of large exponential sums.
% В наших
% последующих
% работах по
% настоящей
% тематике, мы планируем получить дальнейшие
% приложения  результатов о больших тригонометрических суммах к задачам комбинаторной теории чисел
% и, в частности, доказать усиление теоремы \ref{t:Green_3A}.

    The author is grateful to Professor S.V. Konyagin for his helpful idea
    and to N.G. Moshchevitin for constant attention to this work.

\refstepcounter{section}
\label{proofs}

{\bf \arabic{section}. Proof of main result.}

 Let us explain the main idea of the proof of Theorem \ref{t:main}.                         %the?
 Let $N$ be a positive integer and $\F{A}(r)$ the Fourier transform of
 the characteristic function of a set $A$.
 As was noted above, we have Parseval's identity
 \begin{equation}\label{f:2_1}
    \sum_{r\in \Z_N} |\F{A} (r)|^2  = N |A| \,.
 \end{equation}
 Are there another relationships between Fourier coefficients $\F{A}(r)$?
% Существуют ли еще какие--то нетривиальные соотношения между коэффициентами Фурье $\F{A}(r)$
% кроме равенства (\ref{f:2_1})?
 It is easy to see that the answer is yes.

 Let us consider a slightly more general situation.
 Let $f : \Z_N \to \C$ be a complex function.
 We have the inversion formula
 \begin{equation}\label{f:2_2_obr}
    f(x) = \frac{1}{N} \sum_{r\in \Z_N} \F{f} (r) e(rx) \,.
 \end{equation}
 The function $f(x)$ is the characteristic function of a subset of $\Z_N$
 iff for any $x \in \Z_N$ we have
 \begin{equation}\label{f:2_3}
    | f(x) |^2 = f(x) \,.
 \end{equation}
 Combining (\ref{f:2_2_obr}) and  (\ref{f:2_3}), we get
 \begin{equation}\label{}
    \frac{1}{N^2} \sum_{r',r''} \F{f}(r') \overline{\F{f}(r'')} e(r'x - r''x) = \frac{1}{N} \sum_u \F{f}(u) e(ux) \,.
 \end{equation}
 Hence
 \begin{equation}\label{f:2_4}
    \sum_u \left( \frac{1}{N} \sum_r \F{f}(r) \overline{\F{f}(r-u)} \right) e(ux) = \sum_u \F{f}(u) e(ux) \,.
 \end{equation}
 Identity (\ref{f:2_4}) is true for all $x\in \Z_N$.
 It follows that
 \begin{equation}\label{f:2_5}
    \F{f}(u) = \frac{1}{N} \sum_r \F{f}(r) \overline{\F{f}(r-u)} \,.
 \end{equation}
 Thus the complex function $f : \Z_N \to \C$ is the characteristic function
 iff identity (\ref{f:2_5}) holds.
 Moreover, (\ref{f:2_5}) contains all relationships between Fourier coefficients of the set $A$.
 For example to obtain Parseval's identity (\ref{f:Par}) one can put $u=0$
 in formula (\ref{f:2_5}).

 We need in an obvious generalization of (\ref{f:2_5}).
 Suppose that $f,g : \Z_N \to \C$ are two complex functions.
 Then
 \begin{equation}\label{f:2_6}
    \frac{1}{N} \sum_r \F{f}(r) \overline{\F{g}(r-u)} = \sum_x f(x) \overline{g(x)} e(-xu) \,.
 \end{equation}
 Clearly, (\ref{f:2_6}) implies  (\ref{f:2_5}).

 Let us explian the main idea of our proof.
% Приступим к объяснению основной идеи доказательства теоремы \ref{t:main}.
 Let $A\subseteq \Z_N$ be a set, $|A|= \d N$, and
 $\r_\a$ be the set of large exponential sums of $A$.
 Let us consider the model situation.
 Suppose that for all $r\in \r_\a \setminus \{ 0 \}$ we have $|\F{A} (r)| = \a N$,
 and if $r\notin \r_\a$, $r\neq 0$ it follows that $\F{A}(r) = 0$
 (we shall explain this assumption later).
 Let $\d\le 1/4$, and $ 0 \neq u \in \r_\a$ be an arbitrary residual.
 By assumtion $|\F{A} (u)| = \a N$.
 Using (\ref{f:2_5}) and the triangle inequality, we get
 $$
    \a N = |\F{A} (u)| \le \frac{1}{N} \sum_r |\F{A} (r)| |\F{A} (r-u)|
         \le
 $$
 \begin{equation}\label{}
        \le
         \frac{1}{N} \d N |\F{A} (-u)| + \frac{1}{N} |\F{A} (u)| \d N
         + \frac{1}{N} \sum_{r\neq 0, r\neq u} |\F{A} (r)| |\F{A} (r-u)| \,.
 \end{equation}
 It follows that
 $$
    \frac{1}{N} \sum_{r\neq 0, r\neq u} |\F{A} (r)| |\F{A} (r-u)| \ge \frac{\a N}{2} \,.
 $$
 For all $r\neq 0$ we have $|\F{A}(r)| = \a N \cdot \r_\a (r)$.
 Hence
 \begin{equation}\label{f:2_6'}
    \sum_{r\neq 0, r\neq u} \r_\a (r) \r_\a (r-u) \ge \frac{1}{2\a} \,.
 \end{equation}
 Using (\ref{f:2_6'}), we obtain for all $u\in \r_\a \setminus \{ 0 \}$
 the number of solution of the equation $r_1 - r_2 = u$, where $r_1,r_2 \in \r_\a \setminus \{ 0 \}$
 is at least $1/(2 \a)$.
 Therefore  there are  non--trivial arithmetic relationships between the
 elements of the set $\r_\a$.

 Let us give a strict proof of Theorem \ref{t:main}.
 To show the main idea of the proof
 we deduce our result for the simplest case
 $k=2$ and after that in full generality.
 Thus let $k=2$, and
 $B$ be an arbitrary subset of
 $\r_\a \setminus \{ 0 \}$.
 By $[N]$ denote the segment of positive integers $\{1,2,\dots, N\}$.

 We need in the following lemma.

 \Lemma
 \label{l:omega}
 {\it
    Let $\d,\a'$ be real numbers, $0<\a' \le \d$
    and $A$ be a subset of $\Z_N$, $|A| = \d N$.
    Let also
    \begin{equation}\label{f:R'_def}
        \r'_{\a'} = \{~ r\in \Z_N ~:~ \a' N \le |\F{A} (r)| < 2 \a' N ~\}
    \end{equation}
    and $B'$  be an arbitrary subset of $\r'_{\a'} \setminus \{ 0 \}$.
    Then $T_2 (B') \ge (\a')^4  |B'|^4 / (16 \d^3)$.
 }
 \\
 \Proof
 Let
    $$
        f_{B'} (x) = \frac{1}{N} \sum_{r\in B'} \F{A} (r) e(rx) \,.
    $$
 It is easy to see that $\F{f}_{B'} (r) = \F{A} (r) {B'}(r)$.
 Let
 \begin{equation}\label{}
    \sigma = \sum_s | \sum_r \F{f}_{B'} (r) \ov{\F{A} (r-s)} |^2 \,.
 \end{equation}
 Using (\ref{f:2_6}) and Parseval's identity, we get
 \begin{equation}\label{}
    \sigma = N^2 \sum_s | \sum_x f_{B'} (x) \ov{A(x)} e(-xs) |^2 = N^3 \sum_x |f_{B'} (x)|^2 A(x)^2 \,.
 \end{equation}
 Let us obtain a lower bound for $\sum_x |f_{B'} (x)|^2 A(x)^2$.
 Using (\ref{f:Par}) and the definition of the set $\r'_{\a'}$, we have
 \begin{equation}\label{}
    \left( \sum_x f_{B'} (x) A(x) \right)^2 = \left( \frac{1}{N} \sum_r \F{f}_{B'} (r) \ov{\F{A} (r)} \right)^2
    = \left( \frac{1}{N} \sum_r |\F{f}_{B'} (r)|^2 \right)^2
        \ge
 \end{equation}
 \begin{equation}\label{tmp:7_1}
        \ge
            ( N \a'^2 |B'| )^2 = {\a'}^4 |B'|^2 N^2 \,.
 \end{equation}
 On the other hand
 \begin{equation}\label{tmp:7_2}
    \left( \sum_x f_{B'} (x) A(x) \right)^2 \le \left( \sum_x |f_{B'} (x)|^2 A(x)^2 \right)
        \cdot \left( \sum_x A(x)^2 \right) = \d N \left( \sum_x |f_{B'} (x)|^2 A(x)^2 \right).
 \end{equation}
 Using (\ref{tmp:7_1}) and (\ref{tmp:7_2}), we obtain
 \begin{equation}\label{tmp:7_3}
    \sigma^2 \ge \frac{{\a'}^8}{\d^2} |B'|^4 N^8 \,.
 \end{equation}
 Let us obtain an upper bound for $\sigma^2$.
 We have
 $$
    \sigma = \sum_s \sum_{r,r'} \F{f}_{B'} (r) \ov{\F{f}_{B'} (r')}  \cdot \ov{\F{A} (r-s)}  \F{A} (r'-s) =
 $$
 \begin{equation}\label{}
        = \sum_u \left( \sum_r \F{f}_{B'} (r) \ov{\F{f}_{B'} (r-u)} \right) \cdot
            \ov{ \left( \sum_r \F{A} (r) \ov{\F{A} (r-u)} \right) } \,.
 \end{equation}
 It follows that
 \begin{equation}\label{}
    \sigma^2 \le
        \sum_u \left| \sum_r \F{f}_{B'} (r) \ov{\F{f}_{B'} (r-u)} \right|^2
            \cdot
        \sum_u \left| \sum_r \F{A} (r) \ov{\F{A} (r-u)} \right|^2 = \sigma_1 \cdot \sigma_2 \,.
 \end{equation}
 Using (\ref{f:2_5}) and Parseval's identity, we get
 \begin{equation}\label{tmp:7_4}
    \sigma_2 = N^2 \sum_u |\F{A} (u)|^2 = \d N^4 \,.
 \end{equation}
 Since $\F{f}_{B'} (r) = \F{A} (r) {B'}(r)$ and $B' \subseteq \r'_{\a'} \setminus \{ 0\}$
 it follows that
 $|\F{f}_{B'} (r)| \le  2 \a' {B'}(r) N$.
 Hence
 \begin{equation}\label{tmp:7_5}
    \sigma_1 \le 16 (\a')^4 T_2 (B') N^4 \,.
 \end{equation}
 Combining (\ref{tmp:7_4}), (\ref{tmp:7_5}) and (\ref{tmp:7_3}),
 we obtain
 $T_2 (B') \ge (\a')^4  |B'|^4 / (16 \d^3)$.
 This completes the proof of Lemma \ref{l:omega}.
% as required.

 Let
 $$
    B_i = \{~ r\in B ~:~ \a 2^{i-1} N \le |\F{A} (r)| <  \a 2^i N ~\} \,, \quad i \ge 1 \,.
 $$
 Clearly, $B = \bigsqcup_{i\ge 1} B_i$.
 Using Lemma \ref{l:omega}, we obtain
 $T_2 (B_i) \ge (\a 2^{i-1} )^4  |B_i|^4 / (16 \d^3)$, $i\ge 1$.
 Hence
 \begin{equation}\label{t:1_}
    T_2 (B) \ge \sum_i T_2 (B_i) \ge \frac{\a^4}{\d^3 2^8} \sum_i 2^{4i} |B_i|^4 \,.
 \end{equation}
 We have $|B| = \sum_i |B_i|$.
 Using the Cauchy--Schwarz inequality, we get
 \begin{equation}\label{t:2_}
    |B|^4 = \left( \sum_i |B_i| 2^i 2^{-i} \right)^4  \le
        \left( \sum_i 2^{4i} |B_i|^4 \right) \cdot \left( \sum_i 2^{-4i/3} \right)^3
        \le
        %\left(
        \sum_i 2^{4i} |B_i|^4
        %\right)
        \,.
 \end{equation}
 Combining (\ref{t:2_}) and (\ref{t:1_}), we obtain
 \begin{equation}\label{t:3_}
    T_2 (B) \ge \frac{\a^4}{\d^3 2^{8}} |B|^4 \,.
 \end{equation}

 Let us consider the general case $k\ge 2$.

% \Proof
{\bf Proof of Theorem \ref{t:main}}
 First of all let us prove the analog of Lemma \ref{l:omega}.

 \Lemma
 \label{l:omega_k}
 {\it
    Let $\d,\a'$ be real numbers, $0<\a' \le \d$, $A$ be a subset of $\Z_N$, $|A| = \d N$
    and $k\ge 2$ be an even number.
    Let also
    \begin{equation}\label{f:R''_def}
        \r'_{\a'} = \{~ r\in \Z_N ~:~ \a' N \le |\F{A} (r)| < 2 \a' N ~\} \,.
    \end{equation}
    and $B'$ be an arbitrary subset of $\r'_{\a'} \setminus \{ 0 \}$.
    Then $T_k (B') \ge \d (\a')^{2k} |B'|^{2k} / (2 \d)^{2k}$.
 }
 \\
 \Proof
 Let
    $$
        f_{B'} (x) = \frac{1}{N} \sum_{r\in B'} \F{A} (r) e(rx) \,.
    $$
 Consider the sum
 \begin{equation}\label{f:I1}
    \sigma = \left( \sum_x f_{B'} (x) A(x) \right)^k \,.
 \end{equation}
 Clearly,
 \begin{equation}\label{}
    \sigma \ge ( \a'^2 |B'| N )^k
 \end{equation}
 (see proof of Lemma \ref{l:omega}).
 Since $k$ is even it follows that $k=2k'$, $k' \in \N$.
 Using H\"{o}lder's inequality, we get
 $$
    \sigma = \left( \sum_x f_{B'} (x) A(x) \right)^{2k'} \le \left( \sum_x |f_{B'} (x)|^{2k'} A^2(x) \right)
                \left( \sum_x A(x) \right)^{k-1} =
 $$
 \begin{equation}\label{}
                =
                \left( \sum_x |f_{B'} (x)|^{2k'} A^2(x) \right) (\d N)^{k-1} \,.
 \end{equation}
 Hence
 \begin{equation}\label{f:I2}
    \sigma'^2 = \left( \sum_x |f_{B'} (x)|^{2k'} A^2(x) \right)^2 \ge \d^2 \frac{\a'^{4k}}{\d^{2k}} |B'|^{2k} N^2 \,.
 \end{equation}
 On the other hand, using  (\ref{f:2_2_obr}), we have
$$
    \sigma' = \sum_x |f_{B'} (x)|^{2k'} A^2(x) =
        \frac{1}{N^{2k'+2}} \sum_x \sum_{r_1,\dots,r_{k'}, r'_1,\dots,r'_{k'}} \sum_{y,z} \F{f}(r_1) \dots \F{f}(r_{k'})
        \ov{\F{f}(r_1)} \dots \ov{\F{f}(r_{k'})} \F{A} (y) \ov{\F{A}(z)}
$$
$$
        %\cdot
        e(x(r_1+\dots +r_{k'} - r'_1 - \dots -r'_{k'} )) e(x(y-z))
        =
$$
$$
        =
            \frac{1}{N^{2k'+1}} \sum_{u,y}
            \sum_{r_1,\dots,r_{k'}, r'_1,\dots,r'_{k'}, r_1 + \dots + r_{k'} =  r'_1 + \dots + r'_{k'} - u}
            \F{f}(r_1) \dots \F{f}(r_{k'}) \ov{\F{f}(r_1)} \dots \ov{\F{f}(r_{k'})} \F{A} (y) \ov{\F{A}(y-u)}
        =
$$
$$
        =
            \frac{1}{N^{2k'+1}}
            \sum_u \left( \sum_y \F{A} (y) \ov{\F{A}(y-u)} \right) \m
$$
\begin{equation}\label{}
            \m
            \left(
                \sum_{r_1,\dots,r_{k'}, r'_1,\dots,r'_{k'}, r_1 + \dots + r_{k'} =  r'_1 + \dots + r'_{k'} - u}
                \F{f}(r_1) \dots \F{f}(r_{k'}) \ov{\F{f}(r_1)} \dots \ov{\F{f}(r_{k'})}
            \right)
 \end{equation}
 It follows that
 $$
    \sigma'^2  \le \frac{1} {N^{4k'+2}}  \sum_u \left| \sum_y \F{A} (y) \ov{\F{A}(y-u)} \right|^2 \m
 $$
 \begin{equation}\label{}
            \m
            \sum_u
            \left|
                \sum_{r_1,\dots,r_{k'}, r'_1,\dots,r'_{k'}, r_1 + \dots + r_{k'} =  r'_1 + \dots + r'_{k'} - u}
                \F{f}(r_1) \dots \F{f}(r_{k'}) \ov{\F{f}(r_1)} \dots \ov{\F{f}(r_{k'})}
            \right|^2 = \sigma_1 \cdot \sigma_2 \,.
 \end{equation}
 Using (\ref{f:2_5}) and Parseval's identity, we obtain
 \begin{equation}\label{f:I3-}
    \sigma_1 = N^2 \sum_u |\F{A} (u)|^2 = \d N^4 \,.
 \end{equation}
 Since $B' \subseteq \r'_{\a'} \setminus \{ 0 \}$ it follows that $|\F{f}_{B'} (r)| \le 2 \a' B'(r) N$.
 Hence
$$
    \sigma_2 \le \left( (2 \a' N)^{2k'} \right)^2
                    \sum_u
                    \left|
                        \sum_{r_1,\dots,r_{k'}, r'_1,\dots,r'_{k'}, r_1 + \dots + r_{k'} =  r'_1 + \dots + r'_{k'} - u}
                        B'(r_1) \dots B'(r_{k'}) B'(r_1) \dots B'(r_{k'})
                    \right|^2
$$
 \begin{equation}\label{f:I3}
    = (2 \a' N)^{2k} \cdot T_k (B') \,.
 \end{equation}
 Using (\ref{f:I2}), (\ref{f:I3-}) and (\ref{f:I3}), we get
 \begin{equation}\label{}
    T_{k} (B') \ge \d (\a')^{2k} |B'|^{2k} / ( 2 \d)^{2k} \,.
 \end{equation}
 This completes the proof of Lemma \ref{l:omega_k}.

 Let
 % $\l = 2^{1/2k} > 1$ and
 $$
    B_i = \{~ r\in B ~:~ \a 2^{i-1} N \le |\F{A} (r)| <  \a 2^i N ~\} \,, \quad i \ge 1 \,.
 $$
 Clearly, $B = \bigsqcup_{i\ge 1} B_i$.
 Using Lemma \ref{l:omega_k}, we obtain
 $T_k (B_i) \ge \d (\a 2^{i-1} )^{2k}  |B_i|^{2k} / (2 \d)^{2k}$, $i\ge 1$.
 Hence
 \begin{equation}\label{t:1}
    T_k (B) \ge \sum_i T_k (B_i) \ge \frac{\d \a^{2k}}{2^{4k} \d^{2k}} \sum_i 2^{2ki} |B_i|^{2k} \,.
 \end{equation}
 We have $|B| = \sum_i |B_i|$.
 Using H\"{o}lder's inequality, we get
 \begin{equation}\label{t:2}
    |B|^{2k} = \left( \sum_i |B_i| 2^i 2^{-i} \right)^{2k}  \le
        \left( \sum_i 2^{2ki} |B_i|^{2k} \right) \cdot \left( \sum_i 2^{-2k i/(2k-1)} \right)^{2k-1}
        \le
        %\left(
        \sum_i 2^{2ki} |B_i|^{2k}
        %\right)
        \,.
 \end{equation}
 Combining (\ref{t:2}) and
 (\ref{t:1}), we have
 \begin{equation}\label{t:3}
    T_k (B) \ge \frac{\d \a^{2k}}{2^{4k} \d^{2k}} |B|^{2k} \,.
 \end{equation}
 This completes the proof.
% \ref{t:main} доказана.

\refstepcounter{section}
\label{matrix}

\newpage
{\bf \arabic{section}. Linear equations over sets of large exponential sums.}

    Let $k$ be a positive interger, $d\ge 0$ be an integer.
    Let $M = (m_{ij})$ be a matrix $(2^{d+1} k \m (d+1))$, where
    elements $(m_{ij})$ of $M$
    is defined by
\begin{displaymath}
\label{d:aij}
  m_{ij} =
  \left\{ \begin{array}{ll}
    1,            &            \mbox{ if binary expansion of } (j-1) \mbox{ has } 1 \mbox{ on } (i-1) \mbox{ position, } \\
                                    & \mbox{ and } 1\le j \le 2^{d} k \,,
                              \\
    -1,            &            \mbox{ if binary expansion of } (j-1) \mbox{ has } 1 \mbox{ on } (i-1) \mbox{ position, } \\
                                    & \mbox{ and } 2^{d} k < j \le 2^{d+1} k \,,
                              \\
    0,                 &            \mbox{ otherwise. }
  \end{array} \right.
\end{displaymath}
    Recall that binary expansion of a natural number $n$ is defined by
    $n = \sum n_l \cdot 2^{l-1}$, где $l\ge 1$ and $n_l \in \{ 0,1 \}$.

    Let us give an example of $M$. Let $k=2$ and $d=2$. Then
\begin{displaymath}
  M =
  \left( \begin{array}{ccccccccccccccccccc}
    1 & 1 & 1 & 1 & 1 & 1 & 1 & 1 & -1 & -1 & -1 & -1 & -1 & -1 & -1 & -1 \\
    0 & 1 & 0 & 1 & 0 & 1 & 0 & 1 & 0 & -1 & 0 & -1 & 0 & -1 & 0 & -1 \\
    0 & 0 & 1 & 1 & 0 & 0 & 1 & 1 & 0 & 0 & -1 & -1 & 0 & 0 & -1 & -1 \\
  \end{array} \right)
\end{displaymath}

 In this section we prove the following theorem.

\Th
\label{t:matrix}
{\it
    Let $\d,\a$ be real numbers, $0< \a \le \d$,
    $A$ be a subset of $\Z_N$, $|A| = \d N$,
    $k$ be a positive integer, and $d\ge 0$ be an integer.
    Let also $B\subseteq \r_\a \setminus \{ 0 \}$ be an arbitrary set.
    Consider the system of equations
    \begin{equation}\label{f:equations}
        \sum_{i=1}^d \sum_{j=1}^{2^{d+1} k} m_{ij} r_j = 0 \,,
    \end{equation}
    where the elements $(m_{ij})$ of $M = (m_{ij})$
    was defined earlier,
    and $r_j \in B$.
    Then the number of solutions of (\ref{f:equations})
    is at least
    \begin{equation}\label{f:eq_est}
        \left( \frac{\d \a^{2k}}{2^{4k} \d^{2k}} |B|^{2k} \right)^{2^d} \,.
    \end{equation}
}

    Clearly, Theorem \ref{t:matrix} is a generalization of Theorem \ref{t:main}.
    To see this one can put $d=0$ in Theorem \ref{t:matrix}.

    We need in some properties of Gowers uniformity norms in our proof (see \cite{Gow_m}).

 Let $d \ge 0$ be an integer, and
 let
$\{ 0,1 \}^d = \{ \omega = (\omega_1,\dots, \omega_d) ~:~ \omega_j \in \{0,1\}, j=1,2,\dots,d \}$
 be the standard discrete $d$---dimensional cube.
 For any $\omega \in \{0,1\}^d$ we define $|\omega|:= \omega_1 + \dots + \omega_d$.
If $h=(h_1,\dots,h_d) \in \Z_N^d$ we define
$\omega \cdot h := \omega_1 h_1 + \dots \omega_d h_d$.
Let also $\mathcal{C}$ is the conjugation operator
$\mathcal{C} f(x) := \ov{f(x)}$.
By $\| \o \|$ define the sum $\sum_{i=1}^d \o_i \cdot 2^{i-1} + 1$.
Let $\o \in \{ 0,1 \}^d$.
By the same letter $\o$ we define the map $\o : \Z_N^{2^d} \to \Z_N$ by the rule:
if $\v{r} \in \Z_N^{2^d}$ then $\o(\v{r})$ is $\| \o \|$--th component of the vector $\v{r}$.

\Def Let $f: \Z_N \to \C$ be a function.
{\it Gowers uniformity norm} (or Gowers norm)
of $f$ is
\begin{equation}\label{d:Gowers_norm}
    \| f \|_{U^d} :=
                \left( \frac{1}{N^{d+1}}
                \sum_{x\in \Z_N, h\in \Z_N^d} \prod_{\o \in \{ 0,1 \}^d} \mathcal{C}^{|\o|} f (x+\o\cdot h) \right)^{1/2^d} \,.
\end{equation}

We need in the following lemma (see \cite{Gow_m}).

\Lemma {\bf }
\label{l:monoton}
{\it
    Let $f: \Z_N \to \C$ be a function,
    and $d$ be a positive integer.
    Then
    \begin{equation}\label{}
        \| f \|_{U^d} \le \| f \|_{U^{d+1}} \,.
    \end{equation}
}

Let us  prove the analog of Lemma \ref{l:omega_k}.

\Lemma
\label{l:omega_matrix}
{\it
    Let $\d,\a'$ be real numbers, $0<\a' \le \d$,
    $A$ be a subset of $\Z_N$, $|A|= \d N$,
    $k$ be a positive integer, and $d\ge 0$ be an integer.
    Let also
    \begin{equation}\label{f:R'''_def}
        \r'_{\a'} = \{~ r\in \Z_N ~:~ \a' N \le |\F{A} (r)| < 2 \a' N ~\} \,,
    \end{equation}
    and $B'$ be an arbitrary subset of $\r'_{\a'} \setminus \{ 0 \}$.
    Then the number of solutions of system (\ref{f:equations}), where $r_j \in B'$
    is at least
    \begin{equation}\label{f:eq_est_matrix}
        \left( \frac{\d \a'^{2k}}{2^{2k} \d^{2k}} |B'|^{2k} \right)^{2^d} \,.
    \end{equation}
}
\Proof
Let
$$
    f(x) = \frac{1}{N} \sum_{r\in B'} \F{A} (r) e(rx) \,.
$$
   Using H\"{o}lder's inequality, we get
\begin{equation}\label{e:m1}
    \left| \sum_x f(x) A(x) \right|^{2k}
        \le \sum_x |f(x)|^{2k} \cdot \left( \sum_x A(x) \right)^{2k-1}
            = \sum_x |f(x)|^{2k} \cdot (\d N)^{2k-1} \,.
\end{equation}
    On the other hand, using Parseval's identity and the definition of the set $\r'_{\a'}$, we have
\begin{equation}\label{e:m2}
    \sum_x f(x) A(x) = \frac{1}{N} \sum_r \F{f} (r) \ov{\F{A} (r)} = \frac{1}{N} \sum_r |\F{f} (r)|^2 \ge \a'^2 |B'| N \,.
\end{equation}
 Consider the sum
 \begin{equation}\label{f:sigma_begin}
    \sigma = \|\,\, |f|^{2k}\, \|_{U^0} = \|\,\, |f|^{2k}\, \|_{U^1} = \frac{1}{N} \sum_x |f(x)|^{2k} \,.
 \end{equation}
 Combining
 (\ref{e:m1}) and (\ref{e:m2}), we obtain
 \begin{equation}\label{e:sigma_low_bound}
    \sigma \ge \frac{\d \a'^{4k}}{\d^{2k}} |B'|^{2k} \,.
 \end{equation}
 Using Lemma \ref{l:monoton}, we get
 \begin{equation}\label{e:sigma_monoton}
    \sigma^{2^d} \le \frac{1}{N^{d+1}} \sum_{x\in \Z_N}
                        \sum_{h\in \Z_N^d} \prod_{\o \in \{ 0,1 \}^d}  |f (x+\o\cdot h)|^{2k}
                            =
                                \frac{1}{N^{d+1}} \sum_{x\in \Z_N}
                                    \sum_{\v{h}\in \Z_N^d} \left| \prod_{\o \in \{ 0,1 \}^d} f (x+\o\cdot h) \right|^{2k}
 \end{equation}
% Применяя
 By formula (\ref{f:2_2_obr}), we have
 \begin{equation}\label{}
    \prod_{\o \in \{ 0,1 \}^d} f (x+\o\cdot h)
        = \frac{1}{N^{2^d}} \sum_{\v{r} \in \Z_N^{2^d}} \prod_{\o \in \{ 0,1 \}^d} \F{f} (\o(\v{r})) e( \o(\v{r}) (x+\o \cdot h) ) \,.
 \end{equation}
 Hence
 $$
    \sigma^{2^d} =
        \frac{1}{N^{k2^{d+1} + d+1}} \sum_{x\in \Z_N} \sum_{h\in \Z_N^d}
            \sum_{r^{(1)}, \dots, r^{(k)}, r^{(k+1)}, \dots, r^{(2k)} \in \Z_N^{2^d}}
                \prod_{i=1}^{k} \prod_{\o^{(i)} \in \{ 0,1 \}^d} \F{f} (\o^{(i)} (r^{(i)})) e( \o^{(i)} (r^{(i)}) (x+\o^{(i)} \cdot h) )
 $$
 \begin{equation}\label{}
                \m \prod_{i=k+1}^{2k} \prod_{\o^{(i)} \in \{ 0,1 \}^d} \overline{\F{f} (\o^{(i)} (r^{(i)}))} e( - \o^{(i)} (r^{(i)}) (x+\o^{(i)} \cdot h) )
 \end{equation}
 Denote by $\sum$ the following system of linear equations
\begin{eqnarray*}
        \sum_{i=1}^k \sum_{\o^{(i)} \in \{ 0,1 \}^d} \o^{(i)} (r^{(i)}) & = &
                    \sum_{i=k+1}^{2k} \sum_{\o^{(i)} \in \{ 0,1 \}^d} \o^{(i)} (r^{(i)}) \\
        \sum_{i=1}^k \sum_{\o^{(i)} \in \{ 0,1 \}^d,\,\, \o^{(i)}_1 = 1} \o^{(i)} (r^{(i)}) & = &
                    \sum_{i=k+1}^{2k} \sum_{\o^{(i)} \in \{ 0,1 \}^d,\,\, \o^{(i)}_1 = 1} \o^{(i)} (r^{(i)}) \\
        \dots \dots \dots \dots \dots \dots \dots \dots &  & \dots \dots \dots \dots \dots \dots \dots \dots\\
        \dots \dots \dots \dots \dots \dots \dots \dots &  & \dots \dots \dots \dots \dots \dots \dots \dots\\
        \sum_{i=1}^k \sum_{\o^{(i)} \in \{ 0,1 \}^d,\,\, \o^{(i)}_d = 1} \o^{(i)} (r^{(i)}) & = &
                    \sum_{i=k+1}^{2k} \sum_{\o^{(i)} \in \{ 0,1 \}^d,\,\, \o^{(i)}_d = 1} \o^{(i)} (r^{(i)})
\end{eqnarray*}
 We have
$$
    \sigma^{2^d} =
        \frac{1}{N^{k2^{d+1} + d+1}}
            \sum_{r^{(1)}, \dots, r^{(k)}, r^{(k+1)}, \dots, r^{(2k)} \in \Z_N^{2^d}}
                \prod_{i=1}^{k} \prod_{\o^{(i)} \in \{ 0,1 \}^d} \F{f} (\o^{(i)} (r^{(i)}))
$$
$$
                \m \prod_{i=k+1}^{2k} \prod_{\o^{(i)} \in \{ 0,1 \}^d} \overline{\F{f} (\o^{(i)} (r^{(i)}))}
                \sum_{x\in \Z_N} \sum_{h\in \Z_N^d} e( \o^{(i)} (r^{(i)}) (x+\o^{(i)} \cdot h) - \o^{(i)} (r^{(i)}) (x+\o^{(i)} \cdot h) )
$$
\begin{equation}\label{tmp:m_21:49_05}
    = \frac{1}{N^{k2^{d+1}}} \sum_{r^{(1)}, \dots, r^{(k)}, r^{(k+1)}, \dots, r^{(2k)} \in \sum}
            \prod_{i=1}^{k} \prod_{\o^{(i)} \in \{ 0,1 \}^d} \F{f} (\o^{(i)} (r^{(i)}))
            \prod_{i=k+1}^{2k} \prod_{\o^{(i)} \in \{ 0,1 \}^d} \overline{\F{f} (\o^{(i)} (r^{(i)}))}
\end{equation}
 The vectors $r^{(1)}, \dots, r^{(k)}, r^{(k+1)}, \dots, r^{(2k)}$ in (\ref{tmp:m_21:49_05})
 satisfy system $\sum$.
 It is easy to see that this system of equations is system (\ref{f:equations}).
 %, where $r_j \in B'$.

 Indeed, since $\F{f}_{B'} (r) = \F{A} (r) {B'}(r)$ and $B' \subseteq \r'_{\a'} \setminus \{ 0 \}$
 it follows that
 $|\F{f}_{B'} (r)| \le  2 \a' {B'}(r) N$.
 Hence
 \begin{equation}\label{tmp:m_21:58_05}
    \sigma^{2^d}  \le (2^{2k} (\a')^2k )^{2^d} N^{k2^{d+1}} \,.
 \end{equation}
% all components of vectors belongs to B'
 Using (\ref{e:sigma_low_bound}), (\ref{e:sigma_monoton}) and (\ref{tmp:m_21:58_05}), we finally obtain
 \begin{equation}\label{t:B'}
    \sum_{r^{(1)}, \dots, r^{(k)}, r^{(k+1)}, \dots, r^{(2k)} \in \sum}^{*} 1
        \ge \left( \frac{\d (\a')^{4k}}{\d^{2k}} |B'|^{2k} \right)^{2^d} \frac{1}{(2^{2k} (\a')^{2k} )^{2^d}}
        = \left( \frac{\d \a'^{2k}}{2^{2k} \d^{2k}} |B'|^{2k} \right)^{2^d} \,.
 \end{equation}
 All components of all vectors in (\ref{t:B'}) belong to $B'$.
 In other words the number of solutions of system (\ref{f:equations}), where $r_j \in B'$ is at least
 $$
    \left( \frac{\d \a'^{2k}}{2^{2k} \d^{2k}} |B'|^{2k} \right)^{2^d} \,.
 $$
 This completes the proof of Lemma \ref{l:omega_matrix}.

 Let us prove Theorem \ref{t:matrix}.

% \Proof

Let
 $$
    B_i = \{~ r\in B ~:~ \a 2^{i-1} N \le |\F{A} (r)| <  \a 2^i N ~\} \,, \quad i \ge 1 \,.
 $$
 Clearly, $B = \bigsqcup_{i\ge 1} B_i$.

 Let $E$ be a set.
 By $S_{k,d} (E)$ denote the number of solutions of system (\ref{f:equations}), where $r_i \in E$.
 Using Lemma \ref{l:omega_matrix}, we get
 $$
    S_{k,d} (B_i) \ge \left( \frac{\d (\a 2^{i-1})^{2k}}{2^{2k} \d^{2k}} |B_i|^{2k} \right)^{2^d} \,,
 $$
 where $i\ge 1$.
 Hence
 \begin{equation}\label{t:1m}
    S_{k,d} (B) \ge \sum_i S_{k,d} (B_i) \ge
        \left( \frac{\d \a^{2k}}{2^{4k} \d^{2k}} \right)^{2^d} \sum_i \left( 2^{2ki} |B_i|^{2k} \right)^{2^d} \,.
 \end{equation}
 We have $|B| = \sum_i |B_i|$.
 Using H\"{o}lder's inequality, we obtain
 $$
    |B|^{k2^{d+1}} = \left( \sum_i |B_i| 2^i 2^{-i} \right)^{k2^{d+1}}  \le
        \left( \sum_i \left( 2^{2ki} |B_i|^{2k} \right)^{2^d} \right)
            \cdot \left( \sum_i 2^{- (k2^{d+1} i) /(k2^{d+1}-1)} \right)^{k2^{d+1} - 1}
        \le
 $$
\begin{equation}\label{t:2m}
        \le
%        \left(
        \sum_i \left( 2^{2ki} |B_i|^{2k} \right)^{2^d}
%        \right)
        \,.
 \end{equation}
 Combining (\ref{t:2m}) and
 (\ref{t:1m}),
 we have
 \begin{equation}\label{t:3m}
    S_{k,d} (B) \ge \left( \frac{\d \a^{2k}}{2^{4k} \d^{2k}} |B|^{2k} \right)^{2^d} \,.
 \end{equation}
 This completes the proof of Theorem \ref{t:matrix}.

\refstepcounter{section}
\label{applications}

{\bf \arabic{section}. Applications.}

   Chang \cite{Ch_Fr} used in her proof a theorem of Rudin
   \cite{Rudin_book} (see also \cite{Rudin}) concerning
   dissociated subsets of $\Z_N$.
   We say that $\mathcal{D} = \{ d_1, \dots, d_{|\mathcal{D}|} \} \subseteq \Z_N$ is
   {\it dissociated} if the equality
   \begin{equation}\label{}
        \sum_{i=1}^{|\mathcal{D}|} \eps_i d_i = 0 \pmod N \,,
   \end{equation}
   where $\eps_i \in \{ -1,0,1 \}$ implies that all $\eps_i$ are equal to zero.

   \Th {\bf (Rudin)}
   \label{t:Rudin}
   {\it
        There exists an absolute constant  $C>0$ such that
        for any dissociated set $\mathcal{D} \subseteq \Z_N$,
        any complex numbers $a_n \in \C$, and all positive integers $p \ge 2$
        the following inequality holds
        \begin{equation}\label{f:Rudin}
            \frac{1}{N} \sum_{x\in \Z_N} \left| \sum_{n\in \mathcal{D}} a_n e(nx) \right|^p \le
                (C \sqrt{p})^p \left( \sum_{n\in \mathcal{D}} |a_n|^2 \right)^{p/2} \,.
        \end{equation}
   }

   The proofs of Theorem  \ref{t:Rudin} and Chang's theorem can be found in
   \cite{Green_Chang1, Green_Chang2}.
   Let us show that Rudin's theorem and Theorem \ref{t:main}
   imply Theorem \ref{t:Chang}.

   \Pred
   \label{t:Chang_my}
   {\it
        Theorem \ref{t:main} and Theorem \ref{t:Rudin} imply
        Theorem \ref{t:Chang}.
   }
   \\
   \Proof
    Let $k=2 \lceil \log (1/\d) \rceil$,
    and let $\mathcal{D} \subseteq \r_\a$ be a {\it maximal} dissociated subset of $\r_\a$.
    Since $\mathcal{D}$ is a dissociated set it follows that $0\notin \mathcal{D}$.
    Using Theorem \ref{t:main}, we get
    \begin{equation}\label{tmp:3_1}
        T_k (\mathcal{D}) \ge \frac{\d \a^{2k}}{2^{4k} \d^{2k}} |\mathcal{D}|^{2k} \,.
    \end{equation}
    On the other hand
    \begin{equation}\label{tmp:3_2}
        T_k (\mathcal{D}) \le C^{2k} 2^k k^k |\mathcal{D}|^k \,,
    \end{equation}
    where $C$ is an absolute constant (see Theorem \ref{t:Rudin}).
    Indeed, let number $a_n$ in (\ref{f:Rudin}) equals $\mathcal{D} (n)$,
    and $p=2k$.
    Then the left hand side of (\ref{f:Rudin}) is $T_k (\mathcal{D})$, and
    the right hand side is equal to $C^{2k} 2^k k^k |\mathcal{D}|^k$.
    We have $k=2 \lceil \log 1/\d \rceil$.
    Using (\ref{tmp:3_1}) and (\ref{tmp:3_2}), we obtain
    $|\mathcal{D}| \le 2^8 C^2 (\d /\a)^2 (\log 1/\d)$.
    Since $\mathcal{D}$ is a maximal dissociated subset of $\r_\a$ it follows that
    any $r \in \r_\a$ can be written as
    $r=\sum_{i=1}^{|\mathcal{D}|} \eps_i d_i \pmod N$, where
    $d_i \in \mathcal{D}$ and
    $\eps_i \in \{ -1,0,1 \}$.
%    Заметим, что оценка
%    $|\mathcal{D}| \le 2^8 C^2 (\d /\a)^2 \log (1/\d)$
%    отличается от аналогичной оценки в теореме Чанга лишь
%    в константу раз.                                                                                        %пл.
    This completes the proof of Proposition \ref{t:Chang_my}.

   In the paper we obtain an improvement of Chang's result.
   %%Наш метод доказательства имеет много общих черт с
   %%работами \cite{Vinorgadov,Linnik,Nesterenko}.
   We prove our Theorem \ref{t:Chang_log} in the spirit of works
   \cite{Vinorgadov,Linnik,Nesterenko}.

   \Th
   \label{t:Chang_log}
   {\it
    Let $N$ be a positive integer, $(N,6)=1$,
    $\d,\a$ be real numbers, $0< \a \le \d \le 1/16$,
    and
    $A$ be a subset of $\Z_N$, $|A| = \d N$.
    Then there exists a  set
    $\L^* \subseteq \Z_N$,
    \begin{equation}\label{f:L_est}
        |\L^*| \le \min \,( \, \max(\, 2^{30} (\d /\a)^2 \log (1/\d), 2^{ 4 ( \log \log (1/\d) )^2 + 2} \, ),
                                            \,\,  2^{20} (\d/\a)^2 \log^{13/7} (1/\d) \,)
    \end{equation}
    such that for any residual
    $r \in \r_\a$ there exists a tuple $\l_1^*, \dots, \l_M^* \in \L^*$, $M\le 8 \log (1/\d)$
    such that
    \begin{equation}\label{f:r=log}
        r = \sum_{i=1}^M \eps_i \l^*_i \pmod N \,,
    \end{equation}
    where
    %$\l^*_i \in \L^*$ and
    $\eps_i \in \{ -1,0,1 \}$.

    Besides there exists a set $\t{\L} \subseteq \Z_N$,
    \begin{equation}\label{f:tL_est}
        |\t{\L}| \le 2^{20} (\d/\a)^2 \log^{5/3} (1/\d) \log \log (1/\d)
    \end{equation}
    such that for any residual
    $r \in \r_\a$ there exists a tuple $\t{\l}_1, \dots, \t{\l}_M \in \t{\L}$, $M \le 8 \log (1/\d)$
    such that  (\ref{f:r=log}) holds.
   }

   \Cor
   \label{c:Chang_log}
   {\it
        Let $N$ be a positive integer, $(N,6)=1$, and
        $\d,\a$ be real numbers, $0< \a \le \d \cdot 2^{- 2 ( \log \log (1/\d) )^2 }$.
    Then there exists a set
    $\L^* \subseteq \Z_N$, $|\L^*| \le 2^{30} (\d /\a)^2 \log (1/\d)$
    such that for any residual
    $r \in \r_\a$ there exists a tuple $\l_1^*, \dots, \l_M^* \in \L^*$, $M\le 8 \log (1/\d)$
    such that
    $
        r = \sum_{i=1}^M \eps_i \l^*_i \pmod N \,,
    $
    where
    %$\l^*_i \in \L^*$ and
    $\eps_i \in \{ -1,0,1 \}$.
   }

   To prove Theorem \ref{t:Chang_log} we need in some statements and definitions.

   \Def Let $k,s$ be positive integers.
   Consider the family $\L(k,s)$ of subsets of $\Z_N$.
   A set $\L = \{ \l_1, \dots, \l_{|\L|} \}$ belongs to the family $\L(k,s)$
   if the equality
   \begin{equation}\label{}
        \sum_{i=1}^{|\L|} \l_i s_i = 0 \pmod N \,, \quad \l_i \in \L \,, \quad s_i \in \Z \,, \quad
        |s_i| \le s\,, \quad \sum_{i=1}^{|\L|} |s_i| \le 2k \,,
   \end{equation}
   implies that all $s_i$ are equal to zero.

    The definition of $\L(k,1)$ can be found in \cite{Ruzsa_independent}.

    Note that for any $\L \in \L(k,s)$, we have $0 \notin \L$ and $\L \cap -\L = \emptyset$.
    For an arbitrary $\L \in \L(k,s)$, we obtain the following upper bound for $T_k (\L)$.

    \St
    \label{st:sol}
    {\it
        Let $k$, $s$ be positive integers, $s\ge 3$,
        $\L$ be a subset of the family $\L(k,s)$,
        and
        $|\L| \ge k$.
        Then
        \begin{equation}\label{f:T_k_up_estimate}
            T_k (\L)  \le 2^{9k} k^k |\L|^k \cdot 2^{ \frac{2s k (\log k)^2 }{\log (k^{2s} |\L|^{s-2})}} \,.
        \end{equation}
    }

    \Exm Let $\log |\L| \ge \log^2 k$, and $\L$ be an arbitrary subset of the family $\L(k,3)$.
    Using (\ref{f:T_k_up_estimate}), we get
    $T_k (\L) \le 2^{20k} k^k |\L|^k$.
%    It is easy to see that the term $k^k |\L|^k$ is best possible.

    {\bf Proof of Statement \ref{st:sol}. }
    Let $x\in \Z_N$ be a residual.
    By $N_k (x)$ define the number of vectors $(\l_1,\dots, \l_k)$ such that
    all $\l_i$ belong to $\L$ and $\l_1 + \dots + \l_k = x$.
    We have $T_k (\L) = \sum_{x\in \Z_N} N^2_k (x)$.
    Let $s_1,\dots, s_l$ be positive integers such that $s_1 + \dots + s_l = k$.
%    One can Можно считать, для определенности, набор $s_1,\dots, s_l$ упорядоченным
%    по убыванию $s_1 \ge s_2 \ge \dots \ge s_l \ge 1$.
    Let
    $$
        E(s_1,\dots,s_l) (x) = \{ (\l_1,\dots,\l_k) ~:~ \mbox{among } \l_1,\dots,\l_k
        \mbox{ there exist exactly } s_1 \mbox{ residuals equals } \t{\l}_1,
    $$
    $$
        s_2 \mbox{ residuals equals } \t{\l}_2 \,, \dots, s_l \mbox{ residuals equals } \t{\l}_l \mbox{ such that }
        s_1 \t{\l}_1 + \dots + s_l \t{\l}_l = x
        $$
        $$
        \mbox{ and all } \t{\l}_i \mbox{ are different} \}
    $$
    Let us denote the set $E(s_1,\dots,s_l) (x)$ by $E({\v{s}}) (x)$ for simplicity.
    Recall that for $s_1,\dots,s_l$ in the definition of
    $E({\v{s}}) (x) = E(s_1,\dots,s_l) (x)$ the following equality holds $\sum_{i=1}^l s_i = k$.
    We have
    $$
        N_k (x) = \sum_{\v{s}} |E({\v{s}}) (x)| \,.
    $$
    It follows that
    \begin{equation}\label{}
        \sigma = T_k (\L) = \sum_{x\in \Z_N} ( \sum_{\v{s}} |E({\v{s}}) (x)| )^2 \,.
    \end{equation}
    Let $\v{s} = (s_1, \dots, s_l)$,
    and
    $G = G(\v{s}) = \{ i ~:~ s_i \le s \}$,
    $B = B(\v{s}) = \{ i ~:~ s_i > s \}$.
    Then $|G(\v{s})| + |B(\v{s})| = l(\v{s}) = l$.
    We have
    \begin{equation}\label{f:zero}
        l \le k - s |B| \,.
    \end{equation}
    Indeed
    \begin{equation}\label{tmp:5:42}
        k = \sum_{i\in G} s_i + \sum_{i\in B} s_i \ge |G| + (s+1) |B| = l + s |B| \,.
    \end{equation}
    Using (\ref{tmp:5:42}) we obtain (\ref{f:zero}).
    Let
    also
    $$
        l_j = l_j (\v{s}) = |\, \{ i ~:~ s_i = j,~ i\in [l] \} \, | \,, \quad \quad \quad j = 1,2,\dots, r=r(\v{s})  \,.
    $$

    Let us prove the following lemma.

    \Lemma
    \label{l:sol}
    {\it
        1)  For all $\v{s}$, $\sum_{i=1}^l s_i = k$ and for any $x\in \Z_N$, we have
            \begin{equation}\label{}
                |E({\v{s}}) (x)| \le \frac{k!}{s_1! \dots s_l!} |\L|^{|B(\v{s})|} \,.
            \end{equation}
        2)  For all $\v{s}$, $\sum_{i=1}^l s_i = k$ the number of $x\in \Z_N$ such that
            $E({\v{s}}) (x) \neq \emptyset$
            does not exceed $|\L|^l / l_1!$.
    }
    \\
    {\bf Proof of Lemma \ref{l:sol}. }
    $1)~$    Let $(\l_1,\dots, \l_k) \in E({\v{s}}) (x)$ be a tuple.
    Then $\sum_{i=1}^k \l_i = \sum_{i=1}^l s_i \t{\l}_i = x$, where
    $\t{\l}_i \in \{ \l_1,\dots, \l_k \}$ are different.
    Let us consider another tuple $(\l'_1,\dots, \l'_k) \in E({\v{s}}) (x)$,
    $\sum_{i=1}^k \l'_i = \sum_{i=1}^l s_i \t{\l}'_i = x$, where
    $\t{\l}'_i \in \{ \l'_1,\dots, \l'_k \}$ are different.
    Suppose that for all $i\in B(\v{s})$ we have $\t{\l}_i = \t{\l}'_i$.
    This assumption implies that for any $i\in G(\v{s})$
    we have
    $\t{\l}_i = \t{\l}'_i$.
    Indeed $\sum_{i=1}^l s_i \t{\l}_i = x = \sum_{i=1}^l s_i \t{\l}'_i$.
    It follows that $\sum_{i\in G} s_i \t{\l}_i = \sum_{i\in G} s_i \t{\l}'_i$.
    Besides $\L \cap -\L = \emptyset$.
    Hence
    $\sum_{i\in G} s_i \t{\l}_i - \sum_{i\in G} s_i \t{\l}'_i = \sum_{i} s'_i \l^{0}_i = 0$, where
    $s'_i \in \Z$, $|s'_i| \le s$, $\sum_{i} |s'_i| \le 2k$ and $\l^{0}_i \in \L$ are different.
    Using the definition of the family $\L(k,s)$, we obtain that all $s'_i$ are equal to zero.
    It follows that for any $i\in G(\v{s})$
    we have
    $\t{\l}_i = \t{\l}'_i$.
    Hence the tuple $(\l'_1,\dots, \l'_k)$ is a permutation of $(\l_1,\dots, \l_k)$.
    Using the definition of the set $E(\v{s}) (x)$, we obtain that the number
    of these permutations is equal to $k! / (s_1! \dots s_l!)$.
%    Thus, Таким образом, при фиксированных $\t{\l}_i$, $i \in B$ мы имеем не более
%    $k! / (s_1! \dots s_l!)$ наборов $(\l_1,\dots, \l_k)$, принадлежащих $E({\v{s}}) (x)$.
    It follows that the cardinality of $E({\v{s}}) (x)$
    does not exceed $|\L|^{|B(\v{s})|} \cdot k! / (s_1! \dots s_l!)$.

    2)
    For any $\v{s}$,\, $\sum_{i=1}^l s_i = k$ the number of $x\in \Z_N$ such that
    $E({\v{s}}) (x) \neq \emptyset$ does not exceed the number of tuples
    $(\t{\l}_1, \dots, \t{\l}_l)$, $\sum_{i=1}^l s_i \t{\l}_i = x$, where
    $\t{\l}_i$ are different.
    Let $L_1 = \{ i \in [l] ~:~ s_i = 1\}$.
    We have $|L_1| = l_1$.
    Consider all tuples $(\t{\l}_1, \dots, \t{\l}_l)$ and fix all
    $\t{\l}_i$ such that $i\notin L_1$.
    Note that these residuals $\t{\l}_i$ can be chosen in at most  $|\L|^{l-l_1}$ ways.
    Since all $\t{\l}_i$, $i\in L_1$ are different it follows that
    the number of $\t{\l}_i$ does not exceed
    $\binom{|\L|}{l_1} \le |\L|^{l_1} / l_1!$.
    It follows that the number of all tuples
    $(\t{\l}_1, \dots, \t{\l}_l)$, $\sum_{i=1}^l s_i \t{\l}_i = x$,
    where
    $\t{\l}_i$ are different at most $|\L|^{l-l_1} |\L|^{l_1} / l_1! = |\L|^l / l_1!$.
    This completes the proof of Lemma \ref{l:sol}.

%    Вернемся к доказательству утверждения \ref{st:sol}.

    Let $t = (k \log k ) / \log (k^{2s} |\L|^{s-2})$.
    Let us estimate the sum $\sigma$.
    \begin{equation}\label{}
        \sigma \le 2 \left( \sum_{x\in \Z_N} \left( \sum_{\v{s} : |B(\v{s})| \le t} |E({\v{s}}) (x)| \right)^2
            +  \sum_{x\in \Z_N} \left( \sum_{\v{s} : |B(\v{s})| > t} |E({\v{s}}) (x)| \right)^2 \right) = 2 \sigma_1 + 2 \sigma_2 \,.
    \end{equation}
    We have
    \begin{equation}\label{}
        \sigma_2 \le \sum_{\v{s}_1, \v{s}_2 : |B(\v{s}_1)| > t, |B(\v{s}_2)| > t} \sum_x |E({\v{s}_1}) (x)| \cdot  |E({\v{s}_2}) (x)|
    \end{equation}
    If $|B(\v{s}_1)| > |B(\v{s}_2)|$ then put $\v{s}^* = \v{s}_1$.
    If $|B(\v{s}_1)| \le |B(\v{s}_2)|$ then set
    $\v{s}^* = \v{s}_2$.
    Let also $P_k (\v{s}) = k! / (s_1! \dots s_l!)$.
    Using Lemma \ref{l:sol}, we obtain $|E({\v{s}_1}) (x)| \le P_k(\v{s}_1) |\L|^{|B(\v{s}^*)|}$
    and $|E({\v{s}_2}) (x)| \le P_k(\v{s}_2) |\L|^{|B(\v{s}^*)|}$.
    Using the same lemma one more, we get
    \begin{equation}\label{}
       \sigma_2  \le \sum_{\v{s}_1, \v{s}_2 : |B(\v{s}_1)| > t, |B(\v{s}_2)| > t} |\L|^{l(\v{s}^*)} |\L|^{2|B(\v{s}^*)|}
                        P_k(\v{s}_1) P_k(\v{s}_2) \,.
    \end{equation}
    By (\ref{f:zero}), we have
    \begin{equation}\label{}
        \sigma_2 \le \sum_{\v{s}_1, \v{s}_2 : |B(\v{s}_1)| > t, |B(\v{s}_2)| > t}
                        |\L|^{k-s|B(\v{s}^*)|} |\L|^{2|B(\v{s}^*)|} P_k(\v{s}_1) P_k(\v{s}_2)
                            \le
    \end{equation}
    \begin{equation}\label{}
                            \le |\L|^k |\L|^{-t(s-2)} \sum_{\v{s}_1, \v{s}_2} P_k(\v{s}_1) P_k(\v{s}_2)
                                \le 2^4 |\L|^k |\L|^{-t(s-2)} (k^k)^2 \,.
    \end{equation}
    Since $t = (k \log k ) / \log (k^{2s} |\L|^{s-2})$ it follows that
    \begin{equation}\label{f:k^k_L^-t}
        k^k |\L|^{-t(s-2)} \le  2^{ \frac{2s k (\log k)^2 }{\log (k^{2s} |\L|^{s-2})}} \,.
    \end{equation}
    Hence
    \begin{equation}\label{f:A}
        \sigma_2 \le 2^4 k^k |\L|^k 2^{ \frac{2s k (\log k)^2 }{\log (k^{2s} |\L|^{s-2})}} \,.
    \end{equation}
    Let us estimate $\sigma_1$.
    $$
        \sigma_1 \le 2 \left( \sum_{x\in \Z_N} \left( \sum_{\v{s} : |B(\v{s})| \le t, l(\v{s}) \le k - st} |E({\v{s}}) (x)| \right)^2
                                + \sum_{x\in \Z_N} \left( \sum_{\v{s} : |B(\v{s})| \le t, l(\v{s}) > k - st} |E({\v{s}}) (x)| \right)^2 \right)
                                =
    $$
    \begin{equation}\label{}
                                =
                                2\sigma^{'}_1 + 2\sigma^{''}_1 \,.
    \end{equation}
    We have
    \begin{equation}\label{}
        \sigma^{'}_1 \le \sum_{\v{s}_1, \v{s}_2 : |B(\v{s}_1)|, |B(\v{s}_2)| \le t,\,\, l(\v{s}_1), l(\v{s}_2) \le k - st}
                            \sum_x |E({\v{s}_1}) (x)| \cdot  |E({\v{s}_2}) (x)|
    \end{equation}
%    Let обозначение $\v{s}^*$ имеет тот же смысл, что and выше.
    Using Lemma \ref{l:sol} and inequality (\ref{f:k^k_L^-t}), we obtain
    \begin{equation}\label{}
         \sigma^{'}_1 \le \sum_{\v{s}_1, \v{s}_2 : |B(\v{s}_1)|, |B(\v{s}_2)| \le t,\,\, l(\v{s}_1), l(\v{s}_2) \le k - st}
                            |\L|^{l(\v{s}^*)} |\L|^{2|B(\v{s}^*)|} P_k(\v{s}_1) P_k(\v{s}_2)
                            \le
    \end{equation}
    \begin{equation}\label{f:B}
                            \le 2^4 |\L|^k |\L|^{-t(s-2)} (k^k)^2 \le 2^4 k^k |\L|^k
                            2^{ \frac{2s k (\log k)^2 }{\log (k^{2s} |\L|^{s-2})}} \,.
    \end{equation}
    We need in upper bound for $\sigma^{''}_1$.
    For any $\v{s} = (s_1, \dots, s_l)$,\, $\sum_{i=1}^l s_i = k$,
    we have
    \begin{equation}\label{tmp:osel_I}
        l_1 + \dots + l_r = l \quad \mbox{ and } \quad l_1 + 2 l_2 + \dots + r l_r = k \,.
    \end{equation}
    Using (\ref{tmp:osel_I}), we get
    $
        l = k  - ( l_2 + 2l_3 + \dots + (r-1) l_r)
    $.
    On the other hand $l\ge k - st$.
    It follows that $l_2 + 2l_3 + \dots + (r-1) l_r \le st$.
    %Since для всех $\v{s}$ в сумме  $\sigma^{''}_1$ выполнено $l\ge k - st$, то
%    $l_2 + 2l_3 + \dots + (r-1) l_r \le st$.
    Further
    $l_2 + \dots + l_r \le l_2 + 2l_3 + \dots + (r-1) l_r \le st$
%    It follows that
    and
    $l_1 = l - (l_2 + \dots + l_r) \ge l - st \ge k - 2st$.
 %   Let обозначение $\v{s}^*$ имеет тот же смысл, что and выше.
    Using Lemma \ref{l:sol} and (\ref{f:zero}), we get
    \begin{equation}\label{}
        \sigma^{''}_1 \le \sum_{\v{s}_1, \v{s}_2 : |B(\v{s}_1)|, |B(\v{s}_2)| \le t,\, l(\v{s}_1), l(\v{s}_2) > k - st}
                            \sum_x |E({\v{s}_1}) (x)| \cdot  |E({\v{s}_2}) (x)|
                            \le
    \end{equation}
    \begin{equation}\label{}
                            \le
                                \sum_{\v{s}_1, \v{s}_2 : |B(\v{s}_1)|, |B(\v{s}_2)| \le t,\,\, l(\v{s}_1), l(\v{s}_2) > k - st}
                                    \frac{|\L|^{l(\v{s}^*)|} |\L|^{2|B(\v{s}^*)|}}{l_1 (\v{s}^*) !} P_k (\v{s}_1) P_k (\v{s}_2)
                            \le
    \end{equation}
    \begin{equation}\label{}
                            \le
                                \sum_{\v{s}_1, \v{s}_2 : |B(\v{s}_1)|, |B(\v{s}_2)| \le t,\,\, l(\v{s}_1), l(\v{s}_2) > k - st}
                                    \frac{|\L|^{k - s|B(\v{s}^*)|} |\L|^{2|B(\v{s}^*)|}}{l_1 (\v{s}^*) !} P_k (\v{s}_1) P_k (\v{s}_2) \,.
    \end{equation}
    Since $l_1 = l_1 (\v{s}^*) \ge k - 2st$ it follows that
    \begin{equation}\label{}
        \sigma^{''}_1 \le \frac{|\L|^k}{[k-2st]!} \sum_{\v{s}_1, \v{s}_2} P_k (\v{s}_1) P_k (\v{s}_2)
                            \le 2^4 \frac{|\L|^k}{[k-2st]!} (k^k)^2 \,.
    \end{equation}
    By assumption $|\L| \ge k$.
    Hence $t\le k / (3s - 2)$.
    Using the last inequality, we have
    $$
        [k-2st]! \ge [k-2st]^{[k-2st]} / e^{k} \ge k^{[k-2st]} / (8e)^{k} \,.
    $$
    Since $t = (k \log k) / \log (k^{2s} |\L|^{s-2} )$ it follows that
    $k^{2st} \le 2^{ \frac{2s k (\log k)^2}{\log (k^{2s} |\L|^{s-2})}}$.
    Further
    $$
        [k-2st]! \ge k^k / ( 2^{5k} 2^{ \frac{2s k (\log k)^2 }{\log (k^{2s} |\L|^{s-2})}} ) \,.
    $$
    Hence
    \begin{equation}\label{f:C}
        \sigma^{''}_1 \le 2^4 2^{5k} k^k |\L|^k \cdot 2^{ \frac{2s k (\log k)^2 }{\log (k^{2s} |\L|^{s-2})}} \,.
    \end{equation}
    Combining (\ref{f:A}), (\ref{f:B}) and (\ref{f:C}), we finally obtain
    \begin{equation}\label{}
        \sigma = T_k (\L) \le 2^{9k} k^k |\L|^k \cdot 2^{ \frac{2s k (\log k)^2 }{\log (k^{2s} |\L|^{s-2})}} \,.
    \end{equation}
    This completes the proof of Statement \ref{st:sol}.

   {\bf Proof of Theorem \ref{t:Chang_log}}
   Let $k = 2\lceil \log (1/\d) \rceil$.
   Let $\L = \{ \l_1, \dots, \l_{|\L|} \}$ be a maximal subset $\r_\a \setminus \{ 0 \}$
   belongs to the family $\L(k,3)$.
%   Если $\r_\a = \{ 0 \}$, то доказывать нечего.
%   Если $\r_\a \setminus  \{ 0 \}$ не пусто, то тогда and $\L$ не пусто.
   Let $\L^* = (\bigcup_{j=1}^3 j^{-1} \L ) \bigcup \{ 0 \}$.
   Then $|\L^*| \le 4 |\L|$ and $0 \in \L^*$.
   Let us proof that for any $x \in \r_\a \setminus \{ 0 \}$ there exists $j \in [s]$ such that
   \begin{equation}\label{tmp:4:47}
        xj = \sum_{i=1}^{|\L|} \l_i s_i\,, \quad \mbox{ где } s_i \in \Z\,, \quad |s_i| \le s\,, \quad \sum_{i=1}^{|\L|} |s_i| \le 2k \,.
   \end{equation}
   Since for any $i\in [|\L|]$, $j\in [s]$
   we have $j^{-1} \l_i \in \L^*$ it follows that (\ref{tmp:4:47}) implies
   the theorem.
   Thus, let $x$ be an element of $\r_\a \setminus \L$, $x\neq 0$.
   Let us consider all equations $\sum_{i=1}^{|\L|+1} \t{\l}_i s_i = 0$, where
   $\t{\l}_i \in \L \bigsqcup \{ x \}$ and $s_i \in \Z$, $|s_i| \le s$, $\sum_{i=1}^{|\L|+1} |s_i| \le 2k$.
   If all these equations are trivial i.e. we have
   $s_i = 0$, $i\in [|\L|+1]$ then
   we obtain a contradiction with maximality of $\L$.
   It follows that there exists non--trivial equation (\ref{tmp:4:47}) such that not all numbers
   $j, s_1,\dots,s_{|\L|}$ are equal to zero.
   At the same time $j\in [-s, \dots, s]$.
   Since $\L$ belongs to the family $\L(k,3)$ it follows that $j\neq 0$.
   Hence we can assume that $j\in [s]$.
   Since $2k \le 8 \log (1/\d)$ it follows that for any $x\in \r_\a$
   there exists a tuple $\l_1^*, \dots, \l_M^* \in \L^*$, $M\le 8 \log (1/\d)$
   such that (\ref{f:r=log}) holds.

    Let us obtain the bound $|\L^*| \le \max( \, 2^{30} (\d /\a)^2 \log (1/\d), 2^{ 4 ( \log \log (1/\d) )^2 + 2} \, )$.
    %Нам осталось доказать оценку на мощность множества $\L^*$.

    If $\log |\L| < (\log k)^2$ then $|\L| \le 2^{ 4 ( \log \log (1/\d) )^2}$
    and $|\L^*| \le 2^{ 4 ( \log \log (1/\d) )^2 + 2}$.
    Let $\log |\L| \ge (\log k)^2$.
    Using Statement \ref{st:sol}, we get
    $T_k (\L) \le 2^{20k} k^k |\L|^k$.
    On the other hand, using Theorem \ref{t:main}, we obtain
    $T_k (\L) \ge \d \a^{2k} |\L|^{2k} / (2^{4k} \d^{2k})$.
    It follows that $|\L| \le 2^{27} (\d/\a)^2 \log (1/\d)$
    and
    $|\L^*| \le 2^{30} (\d/\a)^2 \log (1/\d)$.

    In any case, we have $|\L^*| \le \max( \, 2^{30} (\d /\a)^2 \log (1/\d), 2^{ 4 ( \log \log (1/\d) )^2 + 2} \, )$.

    Let us prove that $|\L^*| \le 2^{20} (\d/\a)^2 \log^{13/7} (1/\d)$.
    If $|\L| < k$ then $|\L^*| \le 4 |\L| \le 4k$ as required.
    Let $|\L| \ge k$.
    Using Statement \ref{st:sol}, we get
    \begin{equation}\label{}
        T_k (\L) \le 2^{9k} k^k |\L|^k 2^{ \frac{6k \log^2 k}{\log (k^6 |\L|)}}
                        \le 2^{9k} k^k |\L|^k 2^{ \frac{6k \log k}{7} }
                                = 2^{9k} k^k |\L|^k k^{ \frac{6k}{7} } \,.
    \end{equation}
    On the other hand, by Theorem \ref{t:main}, we have
    $T_k (\L) \ge \d \a^{2k} |\L|^{2k} / (2^{4k} \d^{2k})$.
    It follows that $|\L| \le 2^{18} (\d/\a)^2 \log^{13/7} (1/\d)$ and
    $|\L^*| \le 2^{20} (\d/\a)^2 \log^{13/7} (1/\d)$.

    Let us prove the existence of the set $\t{\L}$.
    Let $s = [ \log \log (1/\d) ]$ and $\L_1$ be a maximal subset of $\r_\a \setminus \{ 0 \}$
    such that $\L_1$ belongs to the family $\L(k,s)$, $k = 2\lceil \log (1/\d) \rceil$.
    Let $\t{\L} = \bigcup_{j=1}^s j^{-1} \L_1$.
    We have $|\t{\L}| \le s |\L_1|$.
    It is easy to see that
    for any $r \in \r_\a$ there exists a tuple $\t{\l}_1, \dots, \t{\l}_M \in \t{\L}$, $M\le 8 \log (1/\d)$
    such that (\ref{f:r=log}) holds.

    Finally, let us prove (\ref{f:tL_est}).
    If $|\L_1| < k$ then we are done.
    Let $|\L_1| \ge k$.
    Using Statement \ref{st:sol}, we get
     \begin{equation}\label{}
        T_k (\L) \le 2^{9k} k^k |\L|^k \cdot 2^{ \frac{2s k \log^2 k}{\log (k^{2s} |\L|^{s-2})}}
                        \le 2^{9k} k^k |\L|^k \cdot 2^{ \frac{2sk \log k}{3s-2} }
                                = 2^{9k} k^k |\L|^k k^{ \frac{2sk}{3s-2} } \,.
    \end{equation}
%    С другой стороны из теоремы \ref{t:main} вытекает, что
    On the other hand, by Theorem \ref{t:main}, we have
    $T_k (\L) \ge \d \a^{2k} |\L|^{2k} / (2^{4k} \d^{2k})$.
    It follows that $|\L| \le 2^{20} (\d/\a)^2 \log^{5/3} (1/\d)$ and
    $|\t{\L}| \le 2^{20} (\d/\a)^2 \log^{5/3} (1/\d) \log \log (1/\d)$.
    This completes the proof.

      Let us obtain an application of Theorems \ref{t:main} and \ref{t:Chang_log}.

      Let $K$ be a subset of $\Z_N$ and $\eps \in (0,1)$ be a real number.
      The {\it Bohr set} $B(K,\eps)$ is defined to be the set
      $$
            B(K,\eps) = \{ x\in \Z_N ~:~ \left\| \frac{rx}{N} \right\| < \eps\,, \mbox{ for all }  r\in K \}\,,
      $$
      where $\| \cdot \|$ is the distance to the nearest integer.
      Properties of Bohr sets can be found in J. Bourgain's paper \cite{Bu}.
      In particular Bourgain proved that
      \begin{equation}\label{tmp:Bohr_size_18:37}
        | B(K,\eps) | \ge \frac{1}{2} \eps^{|K|} N \,.
      \end{equation}

      In proof of Freiman's Theorem  \cite{Ch_Fr} (see also \cite{Green_Chang2})
      Chang
      used
      the following proposition.

      \Pred
      \label{p:2A-2A}
      {\it
        Let $N$ be a positive integer, $\d \in (0,1)$ be a real number,
        and $A$ be a subset of $\Z_N$, $|A| = \d N$.
        Then $2A - 2A$ contains some Bohr set $B(K,\eps)$, where
        $|K| \le 8 \d^{-1} \log (1/\d)$ and $\eps = \d / (2^8 \log (1/\d))$.
      }

        We obtain the following improvement of Proposition \ref{p:2A-2A}.

      \Pred
      \label{p:2A-2A'}
      {\it
        Let $N$ be a positive integer, $(N,6)=1$, let $0< \d \le 2^{-256}$ be a real number,
        and let $A$ be a subset of $\Z_N$, $|A| = \d N$.
        Then $2A - 2A$ contains some Bohr set $B(K,\eps)$, where
        $|K| \le 2^{33} \d^{-1} \log (1/\d)$ and $\eps = 1 / (2^8 \log (1/\d))$.
      }

      Using (\ref{tmp:Bohr_size_18:37}) we get the cardinality of the Bohr set
      $B(K,\eps)$ in Proposition \ref{p:2A-2A} is at least
      $(1/2) \cdot 2^{ -8 \d^{-1} (\log 1/\d)^2 } N$.
      In Proposition \ref{p:2A-2A'} the cardinality of the Bohr set is at least
      $(1/2) \cdot 2^{ -2^{35} \d^{-1} (\log 1/\d) (\log \log 1/\d) } N$.

      We need in the following definition.

      \Def Let $f,g : \Z_N \to \C$ be arbitrary functions.
      We define the {\it convolution } $(f*g) (x)$ of the functions $f$ and  $g$ by the formula
      \begin{equation}\label{}
            (f*g) (x) = \sum_{y\in \Z_N} f(y) g(y-x) \,.
      \end{equation}

      It is easy to check
      \begin{equation}\label{f:convolution}
        \F{(f*g)} (r) = \F{f}(r) \ov{\F{g} (r)} \,.
      \end{equation}

      {\bf Proof of Proposition \ref{p:2A-2A'}. }
      Let $\a = \d^{3/2} / 2\sqrt{2}$.
      Using Corollary \ref{c:Chang_log} we find the set
      $\L^* \subseteq \Z_N$, $|\L^*| \le 2^{33} \d^{-1} \log (1/\d)$ such that
      for any $r \in \r_\a$ there exists a tuple  $\l_1^*, \dots, \l_M^* \in \L^*$, $M\le 8 \log (1/\d)$
      such that (\ref{f:r=log}) holds.
      Let $\r_\a^* = \r_\a \setminus \{ 0 \}$.
      Let us consider
      %Bohr set
      $B_1 = B (\r_\a^*, 1 / 20)$.
      For all $x\in B_1$ and any $r\in \r_\a^*$ we have
      \begin{equation}\label{tmp:7:00}
            | 1 - e(rx) | = 2 | \sin(\pi rx/N)| \le \frac{2\pi}{20} < \frac{1}{2} \,.
      \end{equation}
      Clearly, the expression $(A*A*A*A) (x)$ is the
      number of quadruples $(a_1,a_2,a_3,a_4) \in A^4$ such that
      $a_1+a_2 - a_3 - a_4 = x$.
      It follows that, $(A*A*A*A) (x) > 0$ if and  only if  $x\in 2A-2A$.
      Using (\ref{f:2_2_obr}) and (\ref{f:convolution}), we obtain
      $x \in 2A-2A$ if and only if $\sum_r |\F{A}(r)|^4 e(rx) > 0$.
      %We have
      Let $x\in B_1$.
      Then
      $$
        \sum_r |\F{A}(r)|^4 e(rx) = \sum_r |\F{A}(r)|^4 - \sum_r |\F{A}(r)|^4 (1 - e(rx))
        > \frac{1}{2} \sum_r |\F{A}(r)|^4 - 2 \sum_{r\notin R, r\neq 0} |\F{A}(r)|^4 \ge
      $$
      \begin{equation}\label{tmp:7:20}
        \ge
        \frac{1}{2} \d^4 N^4 - 2 \max_{r\notin R, r\neq 0} |\F{A}(r)|^2 \sum_r |\F{A}(r)|^2 \ge
        \frac{1}{2} \d^4 N^4 - 2 \frac{\d^3 N^2}{8} \d N^2 = \frac{\d^4 N^4}{4} > 0 \,.
      \end{equation}
      (we have used Parseval's identity (\ref{f:Par})).
      Inequality (\ref{tmp:7:20}) implies that the set $B_1$ belongs to $2A-2A$.
      Let us consider another Bohr set
      %$B_2 = B(\L^* \setminus \{ 0 \}, 1 / (2^8 \log (1/\d))$
      $B_2 = B(\L^*, 1 / (2^8 \log (1/\d))$
      and let us prove
      $B_2 \subseteq B_1$.
      Since
      for any
      $r \in \r_\a^*$ there exists a tuple $\l_1^*, \dots, \l_M^* \in \L^*$, $M \le 8 \log (1/\d)$
       such that (\ref{f:r=log}) holds
      it follows that
      for all $x\in B_2$, we have
      \begin{equation}\label{}
        \left\| \frac{rx}{N} \right\| \le \sum_{i=1}^M \left\| \frac{\l^*_i x}{N} \right\|
            \le 8 \log (1/\d) \cdot \frac{1}{2^8 \log (1/\d)} < \frac{1}{20} \,.
      \end{equation}
      Thus $B_2 \subseteq B_1$ and we have found  the Bohr set $B_2 \subseteq 2A-2A$
      with the required properties.
      This completes the proof of Proposition \ref{p:2A-2A'}.

\end{document}